\def\th{{^{\rm th}}}
\def\bfe{{\bf e}}
\def\bfv{{\bf v}}
 \def\D{{\partial}}
 \def\BCH{{\rm BCH}}
 
 \def\Aut{{\rm Aut}}
 \def\Diff{{\rm Diff}}
 \def\ad{{\rm ad}}
 \def\half{{\frac{1}{2}}}
\def\third{{\frac{1}{3}}}
\def\twothird{{\frac{2}{3}}}
 \def\ints{{\Bbb Z}}
 \def\lra{{\longrightarrow}}
 \def\takes{{\colon}}
 \def\QED{{\hfill$\Box$}}
 \def\rats{{\bf Q}}
\documentclass{gtart}
\usepackage{enumerate}
\usepackage{latexsym}
\usepackage{amsfonts}
\usepackage{amssymb}
\usepackage{amscd}
\usepackage{amsbsy}
\usepackage{amsmath}
\usepackage{tikz}
\usetikzlibrary{arrows.meta}
\usetikzlibrary{decorations.pathmorphing}

\newcommand*{\Scale}[2][4]{\scalebox{#1}{$#2$}}%

\begin{document}

\title{An explicit symmetric DGLA model of a triangle}

\authors{Itay Griniasty and Ruth Lawrence}

\address{Weizmann Institute and Hebrew University}
\email{itay.griniasty@weizmann.ac.il, ruthel@ma.huji.ac.il}

\begin{abstract}

We give explicit formulae for a differential graded Lie algebra (DGLA) model of the triangle which is symmetric under the geometric symmetries of the cell. This follows the work of Lawrence-Sullivan on the (unique) DGLA model of the interval and of Gadish-Griniasty-Lawrence on an explicit symmetric model of the bi-gon. As in the case of the bi-gon, the essential intermediate step is the construction of a symmetric point. Although in this warped geometry of points given by solutions of the Maurer-Cartan equation and lines given by a gauge transformation by Lie algebra elements of grading zero, the medians of a triangle are not concurrent, various other geometric constructions can be carried out. The construction can similarly be applied to give symmetric model of arbitrary $k$-gons.

 \end{abstract}
 \maketitlepage

 \section {Introduction}
For a regular cell complex $X$, it is possible to associate a DGLA model $A=A(X)$ satisfying the following conditions
 \begin{itemize}
 \item as a Lie algebra, $A(X)$ is freely generated by a set of generators, one for each cell in $X$ and whose grading is one less than the geometric degree of the cell;
 \item vertices (that is $0$-cells) in $X$ give rise to generators $a$ which satisfy the Maurer-Cartan equation $\partial a+{1\over2}[a,a]=0$ (a flatness condition);
 \item for a cell $x$ in $X$, the part of $\D{x}$ without Lie brackets is the geometric boundary $\D_0x$ (where an orientation must be fixed on each cell);
 \item (locality) for a cell $x$ in $X$, $\D{x}$ lies in the Lie algebra generated by the generators of $A(X)$ associated with cells of the closure $\bar{x}$.
 \end{itemize}

 The existence and general construction of such a model was demonstrated by Sullivan in the Appendix to [4].  By [1], there exist consistent (even symmetric) towers of models of simplices, and such towers are unique up to (exact) DGLA isomorphism. The model of an interval is unique [3]. In [2], an explicit symmetric model of the bi-gon (exhibiting the dihedral symmetry of the bi-gon) was given. In this section we collect some general facts about models of cell complexes (see [3]) while in sections 2, 3 we focus on the triangle and its boundary. In section 4 we show how to use a similar procedure for a general $n$-gon.

{\bf General DGLAs.} Recall that a DGLA is a vector space
$A$ over a field $k$ with $\ints$-grading $A=\oplus_{n\in\ints}A_n$ along
with a bilinear map $[.,.]\takes{}A\times{}A\lra{}A$ (bracket, respecting the grading) and a
linear map $\D\takes{}A\lra{}A$ (differential, grading shift $-1$) for which $\D^2=0$
while
 \begin{itemize}
 \item{} (symmetry of bracket) $[b,a]=-(-1)^{|a||b|}[a,b]$;
 \item{} (Jacobi identity) $\ad_{[a,b]}=[\ad_a,\ad_b]$;
 \item{} (Leibniz rule) $[\D,\ad_a]=ad_{\D{a}}$
 \end{itemize}
in terms of the adjoint action of $A$ on itself given by $\ad_e(a)=[e,a]$.

{\bf Points and localisation.} An element $a\in{}A_{-1}$ is called a {\it point} (or said to be {\it flat}) in the model, if it satisfies the Maurer-Cartan equation $\D{a}+{1\over2}[a,a]=0$.  By the {\sl  localisation} of $A$ to a point $a$, denoted $A(a)$, we will mean DGLA which as a graded Lie algebra is
$$\left(\ker\D_a|_{A_0}\right)\oplus\bigoplus_{n>0}A_n$$
with the induced bracket from $A$ and the differential $\D_a\equiv\D+\ad_a$. This contains only non-negative gradings. 

{\bf Edges and flows.} Any element $e\in{}A_0$ defines a {\it flow} on $A$ by
 $$\frac{dx}{dt}=\D{e}-\ad_e(x)\quad{\rm on}\quad A_{-1}\>,\qquad
 \frac{dx}{dt}=-\ad_e(x)\quad{\rm on}\quad A_{\geq0}\>,$$
This flow is called the {\it flow by} $e$, and preserves the grading. In grading $-1$, it preserves flatness, meaning that if the initial condition is at a point ($x(0)$ satisfies Maurer-Cartan) then at all (rational) time its value also satisfies Maurer-Cartan.  Linearity of the differential equation in $e$ ensures that flowing by $e$ for time $t$ is equivalent to flowing by $te$ for a unit time. Denote the result of flowing by $e$ from $a$ for unit time, by $u_e(a)$, so that the solution of the above differential equation satisfies $x(t)=u_{te}(x(0))$.

For a point $a$, the condition that $u_e(a)=a$ is equivalent to $e\in{}A(a)_{0}$ ($e$ is localised at $a$), that is, $\D_ae=0$. This is a linear condition on $e$ and therefore in this case the flow by $e$ fixes $a$ at all time (not only after unit time).

Furthermore (see [2] Lemma 1), if $e$ flows from a point $a$ to a point $b$ in unit time, then $\D_b\circ\exp(-\ad_e)=\exp(-\ad_e)\circ\D_a$ so that $\exp(-\ad_e)$ intertwines the localisation $A(a)$ to the localisation $A(b)$.

{\bf Remark.} The unique DGLA model of an interval has three generators $a$, $b$ of grading $-1$ (from the endpoints) and $e$ of grading $0$ (the 1-cell). The differential is given by the condition $u_e(a)=b$ (see [3]). Explicitly
$$\D{e}=(\ad_e)b+\sum_{i=0}^\infty{\frac{B_i}{i!}}(\ad_e)^i(b-a)={E\over1-e^E}a+{E\over1-e^{-E}}b\>,$$
where $B_i$ denotes the $i\th$ Bernoulli number defined as
coefficients in the expansion
${x\over{e^x-1}}=\sum\limits_{n=0}^\infty{}B_n {x^n\over{}n!}$, $E\equiv\ad_e$ and the expressions in $E$ are considered as formal power series. Similarly, for any 1-cell $e$ in $X$ with endpoints $a$, $b$, it holds that $u_e(a)=b$ in $A=A(X)$.

{\bf BCH.} Denote by $\BCH(x,y)$ (given by the Baker-Campbell-Hausdorff formula)  the unique element in the free Lie algebra on the two generators $x$ and $y$ for which  $(\exp{x}).(\exp{y})=\exp{\BCH(x,y)}$ in the universal enveloping algebra of $A$, or equivalently $$(\exp{\ad_x})\circ(\exp{\ad_y})=\exp{\ad_{\BCH(x,y)}}\in\Aut(A)\>.$$
Then a flow by $e$ for unit time followed by a flow by $f$ for unit time is equivalent to a flow by $\BCH(e,f)$ for unit time,
 $$u_f\circ{}u_e=u_{\BCH(e,f)}\>,$$
at all gradings [3], so that $e\mapsto{}u_e$ is a homomorphism $(A^{(0)},\BCH)\longrightarrow\Diff(A)$.  $\BCH$ is associative  and so a multiple $\BCH$, $\BCH(x_1,\ldots,x_n)$  is also well-defined.

\noindent{\bf Definition}$\>${\sl A piecewise linear path $\gamma$ in $A$, is a sequence of points $a_i$, $0\leq{}i\leq{}m$ in $A$  connected by edges (elements $e_i$, $1\leq{}i\leq{}m$, of $A_0$) which flow between the respective points, $u_{e_i}(a_{i-1})=a_i$ for all $1\leq{}i\leq{}m$.  For such a path, define $\BCH(\gamma)\in{}A_0$ by  $\BCH(\gamma)=\BCH(e_1,\ldots,e_m)$.}

By [1], if $X$ has $c$ connected components and $\{a_1,\ldots,a_c\}$ is a choice of basepoints, one in each connected component, then  the set of points in $A(X)$ is
$$\bigcup_{i=1}^c\big\{u_e(a_i)\bigm|e\in{}A_0\big\}\cup
\big\{u_e(0)\bigm|e\in{}A_0\big\}\>.$$
For each $i$, the map $\pi_i\takes{}e\mapsto{}u_e(a_i)$ is a `fibration', with fibre $\pi_i^{-1}(a_i)$ generated as a vector space by $\{\BCH(\gamma)|\gamma\in\pi_1(X,a_i)\}$, while the map $\pi_0\takes{}e\mapsto{}u_e(0)$ is injective.

\section{The triangle}

Let $\bar\Delta$ be the triangle, with three 0-cells, three 1-cells and one 2-cell. We denote its corresponding model (DGLA) by $\bar{A}$; as a Lie algebra it will be generated freely by $a,b,c$ (grading $-1$), $e,f,g$ (grading 0) and $h$ (grading 1) corresponding to the 0,1,2-cells respectively in $\bar\Delta$.  The differential $\D$ is determined by its values on generators. On vertices, $\D$ is fixed by the Maurer-Cartan condition, e.g. $\D{a}=-\half[a,a]$. On $1$-cells, $\D$ is also unique, e.g
$$\D{e}={E\over1-e^E}b+{E\over1-e^{-E}}c\>,$$
The only freedom in $\bar{A}$ is in $\D{h}\in{}A_0$. The purpose of this paper is to give a formula for $\D{h}$ which is symmetric under the $S_3$ action of the symmetries of the triangle.

Let $\Delta$ denote $\bar\Delta$ with the 2-cell removed, and $A$ its corresponding model, which is unique, $A=\langle{}a,b,c,e,f,g\rangle\subset\bar{A}$.

{\bf Explicit non-symmetric models of the triangle.} The 2-cell with one vertex $\bar{X}^1$, has a model $\bar{A}^1$ with one generator in each degree $-1$,0,1, say $a$, $e$, $h$ respectively with $\D_0e=0$, $\D_0g=e$. The explicit model is
$$\D{}e=[e,a]\>,\qquad\D{}h=e-[a,h]\>.$$
Equivalently, $\D_ae=0$ and $\D_ag=e$. Using the functoriality of the construction $X\mapsto{}A(X)$ under subdivision of intervals, one obtains a model of $\bar\Delta$ (Figure 1) in which
 $$\D{}h=\BCH(g,e,f)-[a,h]\eqno{(1)}$$
This is not symmetric under the symmetries of the triangle (although it is invariant under the reflection in the median from $a$). We could describe this model as `based' at $a$, and will denote it $\bar{A}_a$. Similarly there are models based at the other vertices of the triangle
\begin{align*}
\bar{A}_b:\quad\D{h}&=\BCH(e,f,g)-[b,h]\>,\\
\bar{A}_c:\quad\D{h}&=\BCH(f,g,e)-[c,h]\>.
\end{align*}
The symmetries of the triangle permute $a$, $b$ and $c$. Similarly they permute $e$, $f$, $g$ with an added sign (the sign of the permutation). These symmetries preserve $A$, which was after all the unique model of the triangle boundary $\Delta$. However they permute the three models $\bar{A}_a$, $\bar{A}_b$ and $\bar{A}_c$.

%%%%%%%
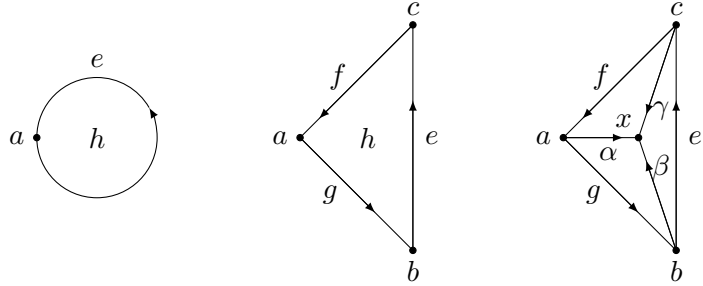
\begin{figure}[h!]
	\centering
	\begin{tikzpicture}
	\begin{scope}[>=latex]	
	\begin{scope}[shift={(-3.5,0)}]
	\draw (0,0) arc(180:30:0.8cm);
	\draw (0,0) arc(-180:30:0.8cm)[->];	
	\draw (0,0) node[anchor=east]{$a$};
	\filldraw (0,0) circle (1.2pt);
	\draw (0.8,0) node {\(h\)};
	\draw (0.8,0.8) node[anchor=south] {\(e\)};
	\end{scope}
	\begin{scope}[shift={(0,0)}]
	\draw (0,0) node[anchor=east]{$a$};
	\draw (1.5,1.5) node[anchor=south]{$c$};
	\draw (1.5,-1.5) node[anchor=north]{$b$};

	\filldraw (0,0) circle (1.2pt);
	\filldraw (1.5,1.5) circle (1.2pt);
	\filldraw (1.5,-1.5) circle (1.2pt);
	
	\draw (0,0) -- (1.5,1.5);
	\draw (1.5,1.5) -- (1/4,1/4)[->];
	\draw (0,0) -- (1.5,-1.5);
	\draw (0,0) -- (1,-1)[->];	
	\draw (1.5,-1.5) -- (1.5,1.5);			
	\draw (1.5,-1.5) -- (1.5,1/2)[->];	
	
	\draw (0.5,0.5) node[anchor=south] {\(f\)};
	\draw (0.4,-0.5) node[anchor=north] {\(g\)};	
	\draw (1.5,0) node[anchor=west] {\(e\)};		
	\draw (1.5/2+0.15,0) node {\(h\)};
	\end{scope}
	
	\begin{scope}[shift={(3.5,0)}]
	\draw (0,0) node[anchor=east]{$a$};
	\draw (1.5,1.5) node[anchor=south]{$c$};
	\draw (1.5,-1.5) node[anchor=north]{$b$};
	
	\filldraw (0,0) circle (1.2pt);
	\filldraw (1.5,1.5) circle (1.2pt);
	\filldraw (1.5,-1.5) circle (1.2pt);
	
	\draw (0,0) -- (1.5,1.5);
	\draw (1.5,1.5) -- (1/4,1/4)[->];
	\draw (0,0) -- (1.5,-1.5);
	\draw (0,0) -- (1,-1)[->];	
	\draw (1.5,-1.5) -- (1.5,1.5);			
	\draw (1.5,-1.5) -- (1.5,1/2)[->];	
	
	\draw (0.5,0.5) node[anchor=south] {\(f\)};
	\draw (0.4,-0.5) node[anchor=north] {\(g\)};	
	\draw (1.5,0) node[anchor=west] {\(e\)};		
	
	\draw (0,0) -- (1,0);	
	\draw (1.5,1.5) -- (1,0);	
	\draw (1.5,-1.5) -- (1,0);
	\draw (0,0) -- (0.8,0)[->];	
	\draw (1.5,1.5) -- (1+0.2*0.5,0.3)[->];	
	\draw (1.5,-1.5) -- (1.1,-0.3)[->];
	
	\filldraw (1,0) circle (1.2pt);
	
	\draw (0.6,0) node[anchor=north] {\(\alpha\)};
	\draw (0.8,0) node[anchor=south] {\(x\)};
	\draw (1.3,0.1) node[anchor=south] 
	{\(\gamma\)};
	\draw (1.3,-0.1) node[anchor=north] 
	{\(\beta\)};
	
	\end{scope}
	\end{scope}
	\end{tikzpicture}
	
	\caption{ {\it Left,}
		The complex \(\bar{X}^1\), can be sub-divided into  \(\bar\Delta\) ({\it center}), however the derived algebra would not be symmetric under \(S_3\).
		{\it Right,}
		the symmetric model \(\bar\Delta\) is based at the central point \(x\), which can be connected to the vertices by \(\alpha,\beta,\gamma\) and are permuted under \(S_3\).	}		
\end{figure}
%%%%%%

{\bf Data for symmetric triangle model.} The aim of this work is to provide an explicit fully symmetric model of $\bar\Delta$. As in [2], this will be done by finding a symmetric point $x$ in $A$ and then producing a model `based' at $x$, meaning that $$\D{h}=q-[x,h]\eqno{(3)}$$ where $q\in\ker\D_x|_{A_0}$ guarantees that $\D^2=0$.  A path from $a$ to $x$ (in $A$) whose $\BCH$ is say $\alpha\in{}A_0$, allows an identification of $A(x)$ with $A(a)$, while  $\ker\D_a|_{A_0}=\rats\BCH(g,e,f)$ generating a model with the correct zeroth order (no Lie brackets) term of $\D$ (namely the topological boundary $\D_0h=e+f+g$) from 
$q=\exp(-\ad_\alpha)\BCH(g,e,f)=\BCH(-\alpha,g,e,f,\alpha)$. The model (3) is symmetric so long as $q$ is anti-symmetric, meaning that it changes by the sign of the permutation under the action of $S_3$.

\section{Construction of symmetric data for triangle}
In this section we work exclusively in the model, $A$, of the triangle boundary $\Delta$ (triangle with 2-cell removed).

{\bf Flattening the triangle.}  For any graph, $\Gamma$, by a {\it realisation} of $\Gamma$ in $A$ we will mean a way of assigning points in $A$ to vertices of $\Gamma$ and elements of $A_0$ to (oriented) edges of $\Gamma$ in such a way that the relation $u_e(a)=b$ holds for every edge of $\Gamma$, where $e$ is assigned to the edge and it connects vertices to which are assigned the points $a$ and $b$. A realisation will be said to be {\it flat}, if the BCH of any loop in the graph in this realisation, vanishes.  

A flat realisation of a connected graph $\Gamma$ is uniquely determined by the label on one vertex, $a$ (an arbitrary point in $A$) and an assignation of elements of $A_0$ to edges in such a way that $BCH(\gamma)=0$ for all loops $\gamma$ in $\Gamma$ based at that vertex (it suffices to check that this holds for a collection of generators $\gamma$ of $\pi_1(\Gamma,a)$). For, given an edge labelling, and the label on one vertex, the realisation condition allows the labels on other vertices to be defined using the flows on edges, and this is always well-defined by the flatness condition.

Thus the graph $\Gamma_0$ with three vertices and three edges has a trivial realisation in $A$, where $a,b,c$ label the vertices and $e,f,g$ label the edges. This is not flat, since there are no relations between $e,f,g$ and in particular $\BCH(e,f,g)\not=0$. 

%%%%%%
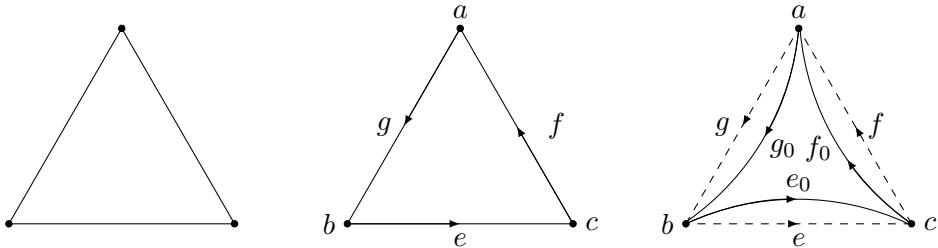
\begin{figure}[h!]
	\centering
	\begin{tikzpicture}
	\begin{scope}[>=latex]
	\begin{scope}[shift={(-4.5,0)}]
	\filldraw (1.5,0) circle (1.2pt);
	\filldraw (-1.5,0) circle (1.2pt);
	\filldraw (0,2.6) circle (1.2pt);
	
	\draw (1.5,0)--(0,2.6);
	\draw (-1.5,0)--(0,2.6);
	\draw (1.5,0)--(-1.5,0);
	\end{scope}
	
	\begin{scope}[shift={(0,0)}]
	\filldraw (1.5,0) circle (1.2pt);
	\filldraw (-1.5,0) circle (1.2pt);
	\filldraw (0,2.6) circle (1.2pt);

	\draw (-1.5,0) node[anchor=east] {\(b\)};
	\draw (1.5,0) node[anchor=west] {\(c\)};	
	\draw (0,2.6) node[anchor=south] {\(a\)};
	
	\draw (1.5,0)--(0,2.6);
	\draw (-1.5,0)--(0,2.6);
	\draw (1.5,0)--(-1.5,0);		
	
	\draw (1.5,0)--(0.75,1.3)[->];
	\draw (0,2.6)--(-0.75,1.3)[->];
	\draw (-1.5,0)--(0,0)[->];

	\draw (0,0) node[anchor=north] {\(e\)};
	\draw (1,1.3) node[anchor=west] {\(f\)};	
	\draw (-0.75,1.3) node[anchor=east] {\(g\)};	
	\end{scope}
	
	\begin{scope}[shift={(4.5,0)}]	
	\filldraw (1.5,0) circle (1.2pt);
	\filldraw (-1.5,0) circle (1.2pt);
	\filldraw (0,2.6) circle (1.2pt);
	
	\draw (-1.5,0) node[anchor=east] {\(b\)};
	\draw (1.5,0) node[anchor=west] {\(c\)};	
	\draw (0,2.6) node[anchor=south] {\(a\)};
	
	\draw [dashed] (0.75+1.5/20,1.3-2.6/20)--(0.75,1.3)[->];
	\draw [dashed] (1.5/20-0.75,1.3+2.6/20)--(-0.75,1.3)[->];
	\draw [dashed] (-1.5/20,0)--(0,0)[->];			

	\draw [dashed] (1.5,0)--(0,2.6);
	\draw [dashed] (-1.5,0)--(0,2.6);
	\draw [dashed] (1.5,0)--(-1.5,0);		
	
	\draw (0,0) node[anchor=north] {\(e\)};
	\draw (0.75,1.3) node[anchor=west] {\(f\)};	
	\draw (-0.75,1.3) node[anchor=east] {\(g\)};
	
	\draw (0,0.3) node[anchor=south] {\(e_0\)};
	\draw (0.6,1) node[anchor=east] {\(f_0\)};	
	\draw (-0.55,1) node[anchor=west] {\(g_0\)};
	
	\draw (0,2.6) arc (-5:-55:3.5cm);
	\draw (0,2.6) arc (-5:-30:3.5cm)[->];
	
	\draw (1.5,0) arc (-90-35:-90-85:3.5cm);
	\draw (1.5,0) arc (-90-35:-90-55:3.5cm)[->];
	
	\draw (-1.5,0) arc (180-65:180-115:3.5cm);
	\draw (-1.5,0) arc (180-65:180-90:3.5cm)[->];		
	\end{scope}
	\end{scope}
	\end{tikzpicture}
	
	\caption{ {\it Left,}
		The triangle graph \(\Gamma_0\) has a trivial (non-flat) realisation in $A$ ({\it center}).
		{\it Right,} A flat realisation of the same graph in $A$.
	}
\end{figure}
%%%%%%

Define $e_0,f_0,g_0\in{}A_0$ by
\begin{align*}
e_0&=\BCH\left(-\third\BCH(e,f,g),e\right)\>,\\
f_0&=\BCH\left(-\third\BCH(f,g,e),f\right)\>,\\
g_0&=\BCH\left(-\third\BCH(g,e,f),g\right)\>.
\end{align*}

{\bf Lemma 3.1}{$\quad{}u_{e_0}(b)=c$}

{\bf Proof} The loop based at $b$ made up of $e,f,g$ in that order verifies that $u_{\BCH(e,f,g)}(b)=b$ and thus  $\BCH(e,f,g)\in\ker\D_b$. Hence also $-\third\BCH(e,f,g)\in\ker\D_b$ so that $u_{-\third\BCH(e,f,g)}(b)=b$. Thus
$u_{e_0}(b)=u_e\left(u_{-\third\BCH(e,f,g)}(b)\right)=u_e(b)=c$, as required.\QED

\vskip 1ex
{\bf Lemma  3.2}{$\quad\BCH(e_0,f_0,g_0)=0$}

{\bf Proof} By properties of $\BCH$, 
$$\exp(-\ad_e)\BCH(e,f,g)=\BCH(-e,\BCH(e,f,g),e)=\BCH(f,g,e)\>.$$
Since $\exp(-ad_e)$ is linear, also $\exp(-\ad_e)\left(-\third\BCH(e,f,g)\right)=-\third\BCH(f,g,e)$ and so $e_0$ can also be written as $e_0=\BCH\left(e,-\third\BCH(f,g,e)\right)$. Thus 
$$\BCH(e_0,f_0)=\BCH\left(e,-\twothird\BCH(f,g,e),f\right)=\BCH\left(e,f,-\twothird\BCH(g,e,f)\right)\>,$$
and combining with $g_0$, $\BCH(e_0,f_0,g_0)=\BCH(e,f,-\BCH(g,e,f),g)=0$.\QED

These two facts together show that we have obtained a flat realisation of $\Gamma_0$ in which $e_0,f_0,g_0$ replace $e,f,g$ on the edges, but the vertices are still assigned $a,b,c$. See Figure 2.

{\bf Iterative step - subdividing a flat triangle.} The graph, $\Gamma_1$, obtained from $\Gamma_0$ by adding midpoints to the edges and joining the three midpoints, wil have six vertices and nine edges. From any flat realisation of $\Gamma_0$, say with edges labelled by $e_0,f_0,g_0$, there can be constructed according to Figure 3, a flat realisation of $\Gamma_1$ in which the corners are labelled by the same points as the given realisation.  To verify flatness, it suffices to verify the condition for the four generating loops around the four smaller triangles in $\Gamma_1$. Verification for the outer triangles is immediate from the definition, while for the inner triangle
$$\BCH\left(\BCH(\half{}f_0,\half{}g_0),\BCH(\half{}g_0,\half{}e_0),\BCH(\half{}e_0,\half{}f_0)\right)=\BCH(\half{}f_0,g_0,e_0,\half{}f_0)$$
which vanishes since $\BCH(e_0,f_0,g_0)=0$, by flatness of the given realisation.

%%%%%%%	
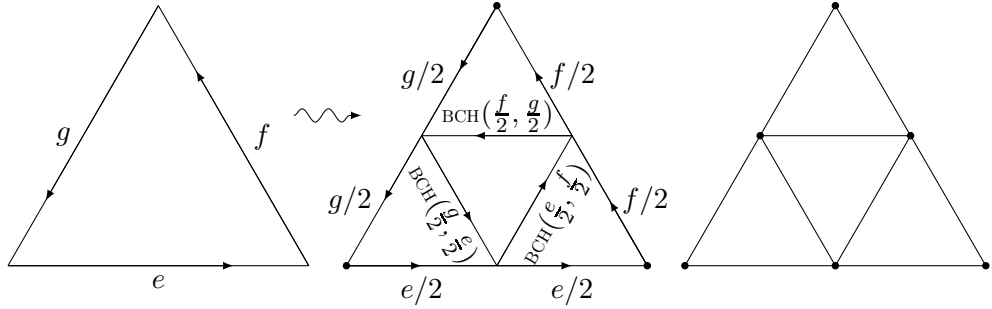
\begin{figure}[h!]
	\centering
	\begin{tikzpicture}
	\begin{scope}[>=latex]
		\draw [decorate,decoration=snake] (-2.7,2) -- (0.9-2.7,2)[->];	
	\begin{scope}[shift={(-4.5,0)}]
		\draw (2,0)--(0,3.46);
		\draw (-2,0)--(0,3.46);
		\draw (2,0)--(-2,0);		
		\draw (2,0)--(0.5,3*1.73/2)[->];	
		\draw (0,3.46)--(-1.5,1.73/2)[->];	
		\draw (-2,0)--(1,0)[->];			
		\draw (0,0) node[anchor=north] {\(e\)};
		\draw (1.1,1.72) node[anchor=west] {\(f\)};		
		\draw (-1,1.72) node[anchor=east] {\(g\)};			
	\end{scope}
	\begin{scope}[shift={(0,0)}]			
		\filldraw (2,0) circle (1.2pt);
		\filldraw (-2,0) circle (1.2pt);
		\filldraw (0,3.46) circle (1.2pt);
		\draw (2,0)--(0,3.46);
		\draw (-2,0)--(0,3.46);
		\draw (2,0)--(-2,0);		
		\draw (2,0)--(1.5,1.73/2)[->];	
		\draw (2,0)--(0.5,3*1.73/2)[->];	
		\draw (0,3.46)--(-1.5,1.73/2)[->];
		\draw (0,3.46)--(-0.5,3*1.73/2)[->];	
		\draw (-2,0)--(-1,0)[->];
		\draw (-2,0)--(1,0)[->];			
		\draw (-1,0) node[anchor=north] {\(e/2\)};
		\draw (1,0) node[anchor=north] {\(e/2\)};
		\draw (1.5,0.87) node[anchor=west] {\(f/2\)};	
		\draw (0.53,2.53) node[anchor=west] {\(f/2\)};	
		\draw (-1.5,0.87) node[anchor=east] {\(g/2\)};	
		\draw (-0.53,2.53) node[anchor=east] {\(g/2\)};	
		\draw (0,0)--(1,1.73);
		\draw (0,0)--(-1,1.73);	
		\draw (-1,1.73)--(1,1.73);
		\draw (1,1.73)--(-0.3,1.73)[->];
		\draw (-1,1.73)--(-1/3,1.73/3)[->];
		\draw (0,0)--(2/3,2*1.73/3)[->];
		\draw (-0.3,-0.1) node[rotate=-60, anchor=east] {\(\Scale[0.6]{\text{BCH}}(\frac{g}{2},\frac{e}{2})\)};			
		\draw (0.3,-0.1) node[rotate=60,anchor=west] {\(\Scale[0.6]{\text{BCH}}(\frac{e}{2},\frac{f}{2})\)};		
		\draw (0,1.65) node[anchor=south] {\(\Scale[0.6]{\text{BCH}}(\frac{f}{2},\frac{g}{2})\)};		
	\end{scope}	
	\begin{scope}[shift={(4.5,0)}]		
		\filldraw (2,0) circle (1.2pt);
		\filldraw (-2,0) circle (1.2pt);
		\filldraw (0,3.46) circle (1.2pt);
		\filldraw (1,1.73) circle (1.2pt);
		\filldraw (-1,1.73) circle (1.2pt);
		\filldraw (0,0) circle (1.2pt);				
		\draw (2,0)--(0,3.46);
		\draw (-2,0)--(0,3.46);
		\draw (2,0)--(-2,0);	
		\draw (0,0)--(1,1.73);
		\draw (0,0)--(-1,1.73);	
		\draw (-1,1.73)--(1,1.73);	
	\end{scope}
	\end{scope}
	\end{tikzpicture}
	\caption{ {\it Left,}
		A flat realisation of \(\Gamma_0\) generating a flat realisation ({\it center}) of \(\Gamma_1\)
		({\it right}).
	}
	
\end{figure}
%%%%%%

{\bf Iterative construction.} Iteratively define $e_n,f_n,g_n$ for non-negative integers $n$, starting with $e_0,f_0,g_0$ defined above, by
$$e_{n+1}=\BCH\left(\half{}f_n,\half{}g_n\right),
f_{n+1}=\BCH\left(\half{}g_n,\half{}e_n\right),
g_{n+1}=\BCH\left(\half{}e_n,\half{}f_n\right)$$
Let $\Gamma_n$ be the graph obtained from $\Gamma_0$ by repeatedly subdividing the inner triangle, $n$ times, each subdivision of the innermost triangle according as the replacement of $\Gamma_0$ by $\Gamma_1$. As in the previous paragraph, starting with a flat realisation of $\Gamma_0$, we obtain a flat realisation of $\Gamma_n$ with the same labels on the corners as the original realisation, and in which the innermost triangle has edges labelled by $e_n,f_n,g_n$. Let $a_n,b_n,c_n$ be the points labelling the vertices of the innermost triangle in $\Gamma_n$. In particular, $a_0=a$, $b_0=b$, $c_0=c$.

Pick any path in $\Gamma_n$ from $a_0$ to $a_n$ and let $\alpha_n\in{}A_0$ denote its BCH in the realisation; this is well-defined since the realisation is flat.

%%%%%%	
\begin{figure}[h!]
	\centering
	\begin{tikzpicture}
	\begin{scope}[>=latex]
	\begin{scope}[shift={(0,0)}]	
	
	\filldraw (3,0) circle (1.2pt);
	\filldraw (-3,0) circle (1.2pt);
	\filldraw (0,5.2) circle (1.2pt);
	
	\filldraw (3/2,5.2/2) circle (1.2pt);
	\filldraw (-3/2,5.2/2) circle (1.2pt);
	\filldraw (0,0) circle (1.2pt);
	
	\filldraw (-3/4,5.2/4) circle (1.2pt);	
	\filldraw (3/4,5.2/4) circle (1.2pt);	
	\filldraw (0,5.2/2) circle (1.2pt);	
	
	\draw (3,0)--(0,5.2);
	\draw (-3,0)--(0,5.2);
	\draw (3,0)--(-3,0);
	
	\draw (0,0)--(3/2,5.2/2);
	\draw (0,0)--(-3/2,5.2/2);	
	\draw (-3/2,5.2/2)--(3/2,5.2/2);		
	
	\draw (-3/4,5.2/4)--(3/4,5.2/4);
	\draw (3/4,5.2/4)--(0,5.2/2);
	\draw (0,5.2/2)--(-3/4,5.2/4);
	
	\draw (3,0) node[anchor=west]{\(c_0\)};
	\draw (-3/2,5.2/2) node[anchor=east]{\(c_1\)};
	\draw (3/4,5.2/4) node[anchor=west]{\(c_2\)};
	
	\draw (-3,0) node[anchor=east]{\(b_0\)};
	\draw (3/2,5.2/2) node[anchor=west]{\(b_1\)};
	\draw (-3/4,5.2/4) node[anchor=east]{\(b_2\)};
	
	\draw (0,5.2) node[anchor=south]{\(a_0\)};
	\draw (0,0) node[anchor=north]{\(a_1\)};
	\draw (0,5.2/2) node[anchor=south]{\(a_2\)};
	
	\draw (3,0)--(3-3/4,5.2/4)[->];	
	\draw (3,0)--(3-3*3/4,3*5.2/4)[->];	
	\draw (0,5.2)--(-3+3/4,5.2/4)[->];	
	\draw (0,5.2)--(-3+3*3/4,3*5.2/4)[->];	
	\draw (-3,0)--(-3/2,0)[->];
	\draw (-3,0)--(3/2,0)[->];
	
	\draw (0,0)--(3/8,5.2/8)[->];
	\draw (0,0)--(3*3/8,3*5.2/8)[->];
	\draw (-3/2,5.2/2)--(-3/2+3/8,5.2/2-5.2/8)[->];	
	\draw (-3/2,5.2/2)--(-3/2+3*3/8,5.2/2-3*5.2/8)[->];	
	\draw (3/2,5.2/2)--(-3/2+3*3/4,5.2/2)[->];	
	\draw (3/2,5.2/2)--(-3/2+3/4,5.2/2)[->];	
	
	\draw (-3/4,5.2/4)--(0,5.2/4)[->];
	\draw (3/4,5.2/4)--(3/8,3*5.2/8)[->];
	\draw (0,5.2/2)--(-3/8,3*5.2/8)[->];						
	
	\draw (-3/2,0) node[anchor=north] {\(e_0/2\)};
	\draw (3/2,0) node[anchor=north] {\(e_0/2\)};
	\draw (3-3/4,5.2/4) node[anchor=west] {\(f_0/2\)};	
	\draw (3-3*3/4,3*5.2/4) node[anchor=west] {\(f_0/2\)};	
	\draw (-3+3*3/4,3*5.2/4) node[anchor=east] {\(g_0/2\)};
	\draw (-3+3/4,5.2/4) node[anchor=east] {\(g_0/2\)};
	
	\draw (-3*3/8+0.1,3*5.2/8-0.3) node[anchor=east] {\(f_1/2\)};
	\draw (-3/8,5.2/8-0.3) node[anchor=east] {\(f_1/2\)};
	\draw (3*3/8-0.1,3*5.2/8-0.3) node[anchor=west] {\(g_1/2\)};
	\draw (3/8-0.1,5.2/8-0.3) node[anchor=west] {\(g_1/2\)};
	\draw (-3/2+3/4,5.2/2) node[anchor=south]{\(e_1/2\)};
	\draw (-3/2+3*3/4,5.2/2) node[anchor=south]{\(e_1/2\)};					
	
	\draw (0,5.2/4) node[anchor=north]{\(e_2\)};
	\draw (3/8,3*5.2/8) node[anchor=west]{\(f_2\)};
	\draw (-3/8,3*5.2/8) node[anchor=east]{\(g_2\)};	
	
	\end{scope}

	\end{scope}
	\end{tikzpicture}
	\caption{ The constructed flat realisation of \(\Gamma_2\)
	}
	
\end{figure}
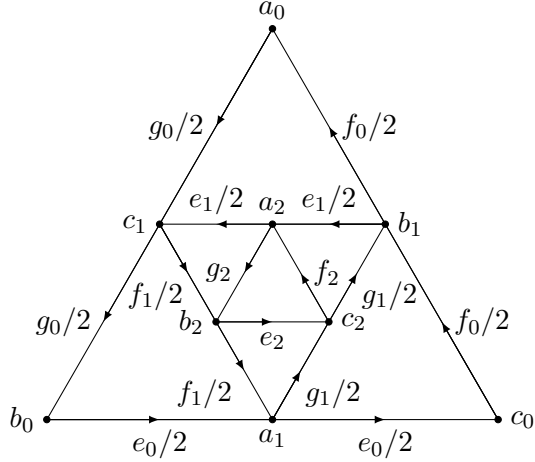
%%%%%%

{\bf Convergence}

{\bf Lemma 3.3}{$\quad{}e_n,f_n,g_n\longrightarrow0$ in $A_0$ as $n\rightarrow\infty$. In other words, for all $m\geq0$,   $e_n^{[m]}\longrightarrow0$ as $n\rightarrow\infty$, where $x^{[m]}\in{}A_0^{[m]}$ denotes the piece of $x\in{}A_0$ with precisely $m$ Lie brackets (and similarly for $f$, $g$).}

{\bf Proof}  Applying the iterative construction above $n$ times to the initial condition $e_1,f_1,g_1$ (in place of $e_0,f_0,g_0$) will arrive at $e_{n+1},f_{n+1},g_{n+1}$. Consequently $e_{n+1},f_{n+1},g_{n+1}$ can be obtained from $e_n,f_n,g_n$ by the replacement $e_0\rightarrow{}e_1$, $f_0\rightarrow{}f_1$, $g_0\rightarrow{}g_1$. Recall that $\BCH(e_0,f_0,g_0)=0$ and so there is a unique Lie algebra expression for $e_n$ as a linear combination of Lie words in $e_0$, $f_0$. Indeed, $e_n,f_n$ can be obtained from $e_0,f_0$ by iterating $n$ times the substitution  
 \begin{align*}
e_0&\longrightarrow{}e_1=\BCH\left(\half{}f_0,\half{}g_0\right)=\BCH\left(\half{}f_0,-\half\BCH(e_0,f_0)\right)\>,\\
f_0&\longrightarrow{}f_1=\BCH\left(\half{}g_0,\half{}e_0\right)=\BCH\left(-\half\BCH(e_0,f_0),\half{}e_0\right)\>.
\end{align*}

Let $B$ be the free Lie algebra on two generators $e_0,f_0$, and consider it embedded in $A_0$ in the natural way. The above substitution induces a linear map $\tau:B\longrightarrow{}B$ which is non-decreasing on the number of Lie brackets and for which $e_n=\tau^n(e_0)$, $f_n=\tau^n(f_0)$. For each $m\geq0$, choose a basis for the finite dimensional vector space $B^{[m]}$. With respect to the basis for $B$ obtained from the union of these bases, the matrix for $\tau$ is a lower triangular (partitioned) matrix.   Since $\tau(e_0)^{[0]}=-\half{}e_0$ and $\tau(f_0)^{[0]}=-\half{}f_0$, thus the diagonal blocks in the matrix of $\tau$ are multiples of the identity matrix with factor $(-2)^{-r}$ on the $r$-th block (dealing with terms with precisely $r-1$ Lie brackets). The truncated (finite-dimensional) matrix of the first $m\times{}m$ blocks gives the matrix of $\tau^{[<m]}$, the induced action of $\tau$ on $B/B^{[\geq{}m]}$. It has eigenvalues $(-2)^{-r}\in(-1,1)$ for $1\leq{}r\leq{}m$,  and thus $\left(\tau^{[<m]}\right)^n\longrightarrow{\bf0}$ as $n\rightarrow\infty$ for all $m$. Applying this to $e_0,f_0$ gives $e_n^{[<m]}\longrightarrow0$ and $f_n^{[<m]}\longrightarrow0$ as $n\rightarrow\infty$; in other words $e_n,f_n\longrightarrow{}0$ in $B$ and hence also in $A_0$. By continuity of $\BCH$, as $g_n=-\BCH(e_n,f_n)$, also $g_n\longrightarrow0$.\QED

\break
{\bf Lemma 3.4}{$\quad$The sequence $(\alpha_n)$ converges in $A_0$.}

{\bf Proof} $\quad$
By the flat realisation of $\Gamma_n$ constructed above, it follows that
$$\alpha_{3n+1}=\BCH\left(\alpha_{3n},\half{}g_{3n},f_{3n+1}\right),\quad 
\alpha_{3n+2}=\BCH\left(\alpha_{3n},\half{}g_{3n},-\half{}e_{3n+1}\right)$$ 
Hence by Lemma 3.3, it suffices to show that the subsequence $(\alpha_{3n})$ converges. Now,
$$\alpha_{3n+3}=\BCH\left(\alpha_{3n},\half{}g_{3n},\half{}f_{3n+1},\half{}e_{3n+2}\right)$$
Let $\sigma$ be the linear map $B\longrightarrow{}B$ defined by the substitution
 \begin{align*}
e_0&\longrightarrow{}g_1=\BCH\left(\half{}e_0,\half{}f_0\right)\>,\\
f_0&\longrightarrow{}e_1=\BCH\left(\half{}f_0,-\half\BCH(e_0,f_0)\right)\>.
\end{align*}
This is the composition of $\tau$ with a rotation. Then $\sigma(e_n)=g_{n+1}$, $\sigma(f_n)=e_{n+1}$ while $\sigma(g_n)=f_{n+1}$ and $\alpha_{3n}=\BCH\left({g_0\over2},\sigma({g_0\over2}),\ldots,\sigma^{3n-1}({g_0\over2})\right)$. Thus it is enough to show that the sequence
$$\left(\BCH\left({g_0\over2},\sigma({g_0\over2}),\ldots,\sigma^{n-1}({g_0\over2})\right)\right)$$
(which contains $\{\alpha_{3n}\}$ as a subsequence) converges, which we do by proving that for any natural number $m$ its projection onto the finite-dimensional vector space $B/B^{[\geq{}m]}$ converges.

\underline{\it Matrix of $\sigma$.}
We use the same notation as in the proof of the previous lemma. The matrix of $\sigma$ is a lower-triangular block matrix. Since $\sigma(e_0)^{[0]}= \half(e_0+f_0)$ and $\sigma(f_0)^{[0]}=-\half{}e_0$, thus the block in the (1,1) position of the partitioned matrix for $\sigma  $ is 
$$\begin{pmatrix}\half&-\half\\ \half&0\end{pmatrix}$$
This is diagonalisable with eigenvalues $-\half\omega,-\half\omega^2$ where $\omega$ is a cube-root of unity.  Choose a basis for $B^{[0]}$ which diagonalises the (1,1) block of $\sigma$ there.

\underline{\it Diagonal blocks of $\sigma$.} Let $\sigma'$ denote the linear map $B\longrightarrow{}B$ induced by the substitution
\begin{align*}
e_0&\longrightarrow{}\half(e_0+f_0)\>,\\
f_0&\longrightarrow{}-\half{}e_0\>.  
\end{align*}
The matrix of $\sigma'$ will be block diagonal and these diagonal blocks will agree with those in $\sigma$. Via the map 
 \begin{align*}
\left(B^{[0]}\right)^{\otimes{}r}&\longrightarrow{}B^{[r-1]}\\
v_1\otimes\cdots\otimes{}v_r&\longmapsto{}[v_1,\ldots,[v_{r-1},v_r]\ldots]
\end{align*}
we can consider $B^{[r-1]}$ as a quotient of $(B^{[0]})^{\otimes{}r}$ by the ideal $I_r$ generated by Jacobi relations. The action of $\sigma'$ on $B^{[0]}$ induces one on $(B^{[0]})^{\otimes{}r}$ preserving $I_r$ and the action on the quotient is precisely the action of $\sigma'$ on $B^{[r-1]}$, described by the $(r,r)$ block in the matrix of $\sigma'$ (or of $\sigma$). By the previous paragraph, $\sigma'|_{B^{[0]}}$ is diagonalisable with eigenvalues $-\half\omega,-\half\omega^2$ and thus the induced action on $(B^{[0]})^{\otimes{}r}$ is also diagonalisable with eigenvalues which all have absolute value $2^{-r}$. The $(r,r)$ block of the matrix for $\sigma$ is a quotient of this and thus also diagonalisable with eigenvalues which all have absolute value $2^{-r}$. 

\underline{\it Bound on matrix entries in powers of $\sigma^{[<m]}$.} 
Fix $m$. We consider only the induced actions on $B/B^{[\geq{}m]}$, that is the first $m\times{}m$ blocks in the matrix representations; let $\sigma^{[<m]}$ denote this induced action from $\sigma$. Choose a basis for $B^{[r-1]}$ which diagonalises the $(r,r)$ block in $\sigma$ for $1\leq{}r\leq{}m$. Let $C$ be the absolute value of the largest matrix entry in $\sigma^{[<m]}$. Let $d_r=\dim{}B^{[r-1]}$ be the size of the $r$-th block.

For any natural number $n$, the matrix for $\sigma^n$ will be a lower triangular block matrix; the diagonal blocks will be diagonal and the entries will have absolute values $2^{-rn}$ in the $(r,r)$ block. The $(a,b)$ entry in the $(i,j)$ block $(i>j)$ of $\sigma^n$ is
$$\sum_{i\geq{}i_1\geq\cdots\geq{}i_{n-1}\geq{}j}\sum_{e_1=1}^{d_{i_1}}\cdots\sum_{e_{n-1}=1}^{d_{i_{n-1}}}
(\sigma_{ii_1})_{ae_1}(\sigma_{i_1i_2})_{e_1e_2}\ldots(\sigma_{i_{n-1}j})_{e_{n-1}b}$$
where $\sigma_{ij}$ denotes the $(i,j)$ block of the partitioned matrix for $\sigma$. For any $i\geq{}i_1\geq\cdots\geq{}i_{n-1}\geq{}j$, let $s_1,\ldots,s_k$ denote the points at which steps occur, that is those $s$ ($1\leq{}s\leq{}n$, in increasing order) for which $i_{s-1}>i_s$ (counting $i_0\equiv{}i$ and $i_n\equiv{}j$). In particular, $i_{s_1-1}=i$ while $i_{s_k}=j$.  The maximum number of steps $k$ is $i-j$. For a particular sequence of steps (that is, where they occur $s_1,\ldots,s_k$ and what are their values $j_1\equiv{}i_{s_1},\ldots,j_{k-1}\equiv{}i_{s_{k-1}}$), the contribution to the above sum is bounded by
$$(2^{-i})^{s_1-1}C(2^{-j_1})^{s_2-s_1-1}C\cdots{}(2^{-j_{k-1}})^{s_k-s_{k-1}-1}C(2^{-j})^{n-s_k}\cdot{}d_{j_1}\cdots{}d_{j_{k-1}}$$
since $\sigma_{ii}$ is diagonal. For fixed $k\leq{}i-j$ and $j_1,\ldots,j_{k-1}$, 
$$\sum_{1\leq{}s_1<\cdots<s_k\leq{}n}(2^{-i})^{s_1-1}(2^{-j_1})^{s_2-s_1-1}\cdots(2^{-j})^{n-s_k} \leq
\begin{pmatrix}n\\k\end{pmatrix}(2^{-j})^{n-k-1}$$
So, if $d=\max\{d_1,\ldots,d_m\}$, an arbitrary entry in the $(i,j)$ block of $\sigma^n$ is bounded by 
$$\sum_{k=1}^{i-j}C^kd^{k-1}\begin{pmatrix}i-j-1\\k-1\end{pmatrix}\begin{pmatrix}n\\k\end{pmatrix}(2^{-j})^{n-k-1}\>.$$
Since $i-j\leq{}m-1$ and $j\geq1$, this bound is at most $2^{-n}$ times a polynomial in $n$ of degree at most $m-1$ and hence all matrix entries in $(\sigma^{[<m]})^n$ can be bounded by $C'(2/3)^n$ for some $C'$ (dependent on $m$).

\underline{\it Bound on coordinates of $v_n\equiv\left(\sigma^n(\half{}g_0)\right)^{[<m]}$.} As above,
$$\left(\sigma^n\left(\half{}g_0\right)\right)^{[<m]}=(\sigma^{[<m]})^n\left(\half{}g_0^{[<m]}\right)$$
which we denote by $v_n\in{}B^{[<m]}$. The matrix elements in the power of $\sigma^{[<m]}$ are all bounded by a multiple of $(2/3)^n$ while the vector $g_0^{[<m]}$ is constant. Thus, in any chosen basis for $B^{[<m]}$,  $v_n$ has all coordinates (and thus also their sum) bounded by a constant (dependent on $m$) times $(2/3)^n$. 

\underline{\it Coefficients in BCH.} From now onwards we will revert to a basis for $B^{[r]}$ in which the basis elements are (a subset of) Lie monomials in $e_0,f_0$ with $r$ brackets. The formula for $\BCH(x,y)$ is an element of the free Lie algebra on $x$ and $y$.
Since $g_0=-\BCH(e_0,f_0)$, the coefficients in the formula are given precisely by the coordinates of $-g_0$ with repsect to the chosen basis. Denote these coefficients $h^{[r]}_j\in\rats$, so that
$$\BCH(e_0,f_0)=\sum_{r=0}^\infty\sum_{j=1}^{d^{[r]}}h^{[r]}_j\bfe^{[r]}_j$$
where $\bfe^{[r]}_j$ is the $j$-th basis vector in $B^{[r]}$ and $d^{[r]}\equiv{}d_{r+1}$ is the dimension of $B^{[r]}$. For example, $d^{[0]}=2$, take $\bfe^{[0]}_1=e_0$, $\bfe^{[0]}_2=f_0$ as basis for $B^{[0]}$, and then $h^{[0]}_1=h^{[0]}_2=1$. Similarly $d^{[1]}=1$, $\bfe^{[1]}_1=[e_0,f_0]$ and $h^{[1]}_1=\half$. For second order brackets, $d^{[2]}=2$, use $\bfe^{[2]}_1=[e_0,[e_0,f_0]]$, $\bfe^{[2]}_2=[f_0,[e_0.f_0]]$ and then $h^{[2]}_1=-h^{[2]}_2={1\over12}$.

\underline{\it Bound on growth of $\BCH$.} Since $\BCH$ is non-decreasing on the number of Lie brackets, it induces a well-defined (associative) binary operation on $B/B^{[\geq{}m]}$. Define a metric on $B/B^{[\geq{m}]}$ by
$$\left|\left|\sum_{r=0}^\infty\sum_{j=1}^{d^{[r]}}a^{[r]}_j\bfe^{[r]}_j\right|\right|=
\sum_{r=0}^{m-1}\sum_{j=1}^{d^{[r]}}|a^{[r]}_j|$$
Let $D$ denote the maximum norm of all Lie monomials in $e_0,f_0$ with at most $m-1$ brackets. For $a\in{}B$, denote by $a^{[r]}\in{}B^{[r]}$ the part of $a$ with $r$ Lie brackets. Then for any $a,b\in{}B$,
$$\left(\BCH(a,b)\right)^{[r]}=\sum_{i=0}^r\sum_{j=1}^{d^{[i]}}h^{[i]}_j\left(\bfe^{[i]}_j(a,b)\right)^{[r]}$$
where $\bfe(a,b)$ is the result of substituting $a,b$ in place of $e_0,f_0$ in the Lie monomial $\bfe\in{}B$. For example
\begin{align*}
\BCH(a,b)^{[0]}&=a^{[0]}+b^{[0]}\>,\\
\BCH(a,b)^{[1]}&=a^{[1]}+b^{[1]}+\half[a^{[0]},b^{[0]}]\>,\\
\BCH(a,b)^{[2]}&=a^{[2]}+b^{[2]}+\half[a^{[0]},b^{[1]}]+\half[a^{[1]},b^{[0]}]+{1\over12}[a^{[0]},[a^{[0]},b^{[0]}]]
-{1\over12}[b^{[0]},[a^{[0]},b^{[0]}]]\>.
\end{align*}
But for any monomial $\bfe\in{}B$ involving $k$ times $e_0$ and $l$ times $f_0$ ($k,l>0$),
$$||\bfe(a,b)^{[<m]}||\leq{}D||a||^k||b||^l\>,$$
since substituting monomials for $e_0,f_0$ in a monomial will produce another monomial, which will have norm at most $D$.
Thus there exist homogeneous polynomials $p_r$ in two variables, of degree $r+1$, such that for all $a,b$,
$$||\BCH(a,b)-a-b||\leq\sum_{r=1}^{m-1}{}p_r(||a||,||b||)$$
with $p_1(x,y)=Dxy/2$, $p_2(x,y)=Dxy(x+y)/12$ and furthermore $p_r(x,y)$ is divisible by $xy$ for all $r$. So in particular,
$$||\BCH(a,b)-a||\leq||b||Q(||a||,||b||)$$ for a suitable polynomial $Q$ in two variables of degree $m-1$.

\underline{\it $\BCH$-Cauchy.} By the previous paragraphs, we have a sequence of vectors $v_n\in{}B/B^{[\geq{}m]}$ satisfying $||v_n||\leq{}D(2/3)^n$ for all $n$ (some constant $D$) and the proof of the lemma will be complete once it is shown that the sequence
$$\left(\BCH(v_0,v_1,\ldots,v_{n-1})^{[<m]}\right)$$
converges in $B^{[<m]}$. Let $X$ denote the maximum value of $Q(x,y)$ when $0\leq{}x,y\leq{}D$. By the previous paragraph, $$||\BCH(a,b)-a||\leq{}X||b|| \quad{\rm whenever}\quad ||a||,||b||\leq{}D\>.$$ Choose $N$ sufficiently large that $(3/2)^N\geq1+2X$. For arbitrary $m\geq{}N$ we see inductively that for any $i\geq{}m$,
$$||\BCH(v_m,\ldots,v_i)||\leq{}D$$
For, this holds when $i=m$ as $||\bfv_m||\leq{}D(2/3)^m\leq{}D$. Assuming it holds for all $i<n$, then 
$$||\BCH(v_m,\ldots,v_i,v_{i+1})-\BCH(v_m,\ldots,v_i)||\leq{}X||v_{i+1}||\leq{}DX(2/3)^{i+1}$$
Combining with the triangle inequality for $i=m,m+1,\ldots,n-1$,
$$||\BCH(v_m,\ldots,v_n)||\leq||v_m||+\sum_{i=m}^{n-1}DX(2/3)^{i+1}\leq{}D(2/3)^m(1+2X)\leq{}D(2/3)^{m-N}$$ 
which is at most $D$, proving the inductive step. Furthermore, for any $n\geq{}m\geq{}N$, 
\begin{align*}
||&\BCH(v_N,\ldots,v_n)-\BCH(v_N,\ldots,v_m)||\\
&=||\BCH(\BCH(v_N,\ldots,v_m),\BCH(v_{m+1},\ldots,v_n))-\BCH(v_N,\ldots,v_m)||\\
&\leq{}X||\BCH(v_{m+1},\ldots,v_n)||\leq{}XD(2/3)^{m+1-N}
\end{align*}
and therefore the sequence $\left(\BCH(v_N,\ldots,v_n)^{[<m]}\right)$ is a Cauchy sequence in $B^{[<m]}$ and hence converges. Since $\BCH^{[<m]}$ is continuous, taking a $\BCH$ of the sequence with $\BCH(v_0,\ldots,v_{N-1})$, will produce a convergent sequence also, namely $\left(\BCH(v_0,\ldots,v_n)^{[<m]}\right)$, as required. 
\QED

Denote the limit of the sequence $\{\alpha_n\}$ by $\alpha$. Set $x=u_\alpha(a)$. Since $a_n=u_{\alpha_n}(a)$, thus $a_n\longrightarrow{}x$.

{\bf Symmetry.} The symmetry group $S_3$ of the triangle permutes the vertices $a,b,c$ and the edges (with signs) $e,f,g$. By construction, $e_0,f_0,g_0$ will be identically permuted (with signs) as $e,f,g$ and the symmetry of the iterative step guarantees that this holds also for $e_n,f_n,g_n$ for all $n$ and finally that $a_n,b_n,c_n$ will be permuted amongst themselves, and similarly for $\alpha_n,\beta_n,\gamma_n$.  

Since $b_n=u_{g_n}(a_n)$,  $g_n\rightarrow0$ (Lemma 3.3) and $a_n\rightarrow{}x$ (Lemma 3.4), thus  $b_n\longrightarrow{}x$. Similarly $c_n\longrightarrow{}x$. Since $S_3$ permutes $a_n,b_n,c_n$, thus $x$ is invariant under this action, that is, it is a symmetric point.

As in Lemma 3.4, also $\{\beta_n\}$ and $\{\gamma_n\}$ are convergent sequences; denote their limits by $\beta,\gamma\in{}A_0$. Since $S_3$ permutes $\alpha_n,\beta_n,\gamma_n$, thus it also permutes $\alpha,\beta,\gamma$. Finally, applying $S_3$ to $u_\alpha(a)=x$ we obtain that also $u_\beta(b)=u_\gamma(c)=x$. 

%%%%%%	
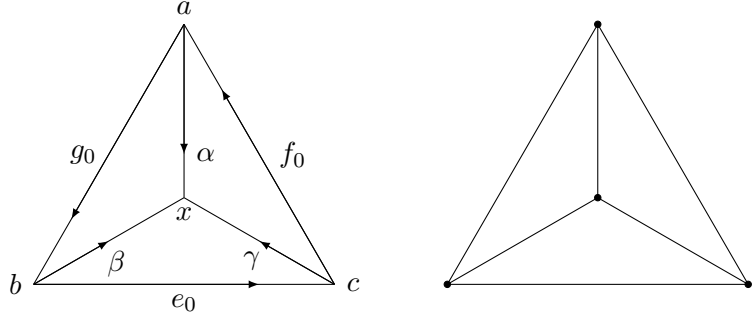
\begin{figure}[h!]
	\centering
	\begin{tikzpicture}
	\begin{scope}[>=latex]
	\begin{scope}[shift={(0,0)}]	
%	\begin{scope}[shift={(7.5,1)}]	
	
	\draw (2,0)--(0,3.46);
	\draw (-2,0)--(0,3.46);
	\draw (2,0)--(-2,0);
	
	\draw (2,0)--(0.5,3*1.73/2)[->];	
	\draw (0,3.46)--(-1.5,1.73/2)[->];	
	\draw (-2,0)--(1,0)[->];			
	
	\draw (0,3.46) node [anchor=south]{\(a\)};		
	\draw (-2,0) node [anchor=east]{\(b\)};		
	\draw (2,0) node [anchor=west]{\(c\)};
	
	\draw (0,0) node[anchor=north] {\(e_0\)};
	\draw (1.1,1.72) node[anchor=west] {\(f_0\)};		
	\draw (-1,1.72) node[anchor=east] {\(g_0\)};
	
	\draw (0,3.46)--(0,3.46/3);
	\draw (-2,0)--(0,3.46/3);
	\draw (2,0)--(0,3.46/3);
	
	\draw (0,3.46)--(0,3.46/2)[->];
	\draw (-2,0)--(-1,3.46/6)[->];
	\draw (2,0)--(1,3.46/6)[->];
	
	\draw (0,3.46/3) node[anchor=north]{\(x\)};	
	\draw (0,3.46/2) node[anchor=west]{\(\alpha\)};	
	\draw (-0.9,3.46/6) node[anchor=north]{\(\beta\)};
	\draw (0.9,3.46/6) node[anchor=north]{\(\gamma\)};
	\end{scope}
	
		\begin{scope}[shift={(5.5,0)}]			
	\draw (2,0)--(0,3.46);
	\draw (-2,0)--(0,3.46);
	\draw (2,0)--(-2,0);
	
	\filldraw (2,0) circle (1.2pt);
	\filldraw (-2,0) circle (1.2pt);
	\filldraw (0,3.46) circle (1.2pt);
	\filldraw (0,3.46/3) circle (1.2pt);
	
	\draw (0,3.46)--(0,3.46/3);
	\draw (-2,0)--(0,3.46/3);
	\draw (2,0)--(0,3.46/3);
	\end{scope}		
	\end{scope}
	\end{tikzpicture}
	\caption{ {\it Left,}
		The constructed flat realisation of the tetrahedron \(T\)
		({\it right}).
	}
	
\end{figure}
%%%%%%

{\bf Symmetric model of the triangle.}
The conclusion from the previous subsections is that $x,\alpha,\beta,\gamma$ constitute symmetric data from which a symmetric realisation of the tetrahedral graph $T$ in $A$ is obtained, with corners $a,b,c$ and internal vertex $x$. The external edges are labeled $e_0,f_0,g_0$ and the internal edges $\alpha,\beta,\gamma$. To see that this is a flat realisation, it must be checked that the BCHs of each of the three generating loops vanish. Note that $\BCH(g_0,\beta_n,-g_n,-\alpha_n)=0$ since it is represented by a loop based at $a$ on the flat realisation of $\Gamma_n$ constructed above. In the limit $n\rightarrow\infty$, the equality gives $\BCH(g_0,\beta,-\alpha)=0$. Similarly for the other faces.

As in the discussion at the end of \S{2}, set $q=\BCH(-\alpha,g,e,f,\alpha)$. 

{\bf Theorem}{$\quad\D{h}=\BCH(-\alpha,g,e,f,\alpha)-[x,h]$ defines a symmetric model of $\bar\Delta$.}
{\bf Proof} $\quad$ It is already known that $x$ is a symmetric point and so it remains to prove (see end of \S{2}) that $q$ is anti-symmetric under the $S_3$ action, for which it is enough to check the action under generators of $S_3$. 

Reflection in the median through $a$ acts by fixing $a$, interchanging $b,c$, changing the sign of $e$, interchanging $f,-g$. This fixes $\alpha$ and interchanges $\beta,\gamma$. This reverses the sign of $\BCH(g,e,f)$ and thus also of $q$.

Rotation cycles between $a,b,c$ and similarly $e,f,g$, $\alpha,\beta,\gamma$. Thus $q$ transforms to 
$\BCH(-\beta,e,f,g,\beta)$. Since $\beta=\BCH(-g_0,\alpha)$, thus
$$\BCH(-\beta,e,f,g,\beta)=\BCH(-\alpha,g_0,e,f,g,-g_0,\alpha)=q$$ where the last step follows, using the definition of $g_0$, from 
$$\BCH(g_0,e,f,g,-g_0)=\BCH(-\third\BCH(g,e,f),g,e,f,\third\BCH(g,e,f))=\BCH(g,e,f)$$
\QED

\section{Generalisations}

{\bf Computations.} By iteratively solving the condition $\sigma(\alpha)=\BCH(-g_0/2,\alpha)$ along with the requirement that $\beta$ is obtained from $\alpha$ (and $\gamma$ from $\beta$) under the rotation $e_0\longrightarrow{}f_0$, $f_0\longrightarrow-\BCH(e_0,f_0)$,
one can calculate $\alpha$, $\beta$, $\gamma$ in terms of $e_0$, $f_0$. The result is
\begin{align*}
\alpha=&-\third(e_0+2f_0)-{1\over6}[e_0,f_0]-{1\over54}[e_0,[e_0,f_0]]+{1\over36}[f_0,[e_0,f_0]]+\cdots\>,\\
\beta=&\third(2e_0+f_0)+{1\over6}[e_0,f_0]+{1\over36}[e_0,[e_0,f_0]]-{1\over54}[f_0,[e_0,f_0]]+\cdots\>,\\
\gamma=&\third(f_0-e_0)-{1\over108}[e_0+f_0,[e_0,f_0]]+\cdots\>.
\end{align*}

{\bf Remark} Note that $\alpha,\beta$ freely generate $B=\langle{}e_0,f_0\rangle$ and so $\gamma$ can be written as a universal Lie word in $\alpha,\beta$, say $\gamma=f(\alpha,\beta)$. The symmetry constraints imply that 
$f(\beta,\alpha)=f(\alpha,\beta)$ while $f(\alpha,f(\alpha,\beta))=\beta$. In fact $f(\alpha,\beta)=-\alpha-\beta+\cdots$ where the first non-trivial term has at four Lie brackets:
$$\frac{17}{2^2\cdot3^3\cdot5\cdot11}\left(A^4\beta+B^4\alpha-A^2B^2\alpha- B^2A^2\beta+\frac{1}{2}(AB^3\alpha+BA^3\beta)\right)\>.
$$
Here $A\equiv\ad_\alpha$ and $B\equiv\ad_\beta$.

{\bf $k$-gons} The arguments of this paper can be applied to any $k$-gon, where the iterative operation is to replace a $k$-gon by inscribing another $k$-gon joining the edge midpoints.  The only slight complication is in the convergence argument. For example, for a square, $\tau$ is replaced by an automorphism of the free Lie algebra on three generators given by 
$$e\rightarrow\BCH\left({e\over2},{f\over2}\right),\quad
f\rightarrow\BCH\left({f\over2},{g\over2}\right),\quad
g\rightarrow\BCH\left({g\over2},-\half\BCH(e,f,g)\right)$$
To zeroth order, this is $e\longmapsto\half(e+f)$, $f\longmapsto\half(f+g)$, $g\longmapsto-\half(e+f)$ which has eigenvalues $0,\half(-1\pm{}i)$ which still all have absolute value less than 1.

\subsection*{ Acknowledgments}
Itay Griniasty is grateful to the Azrieli Foundation for the award of an Azrieli Fellowship. This research was supported by Grant No 2016219 from the United States-Israel Binational Science Foundation (BSF).

\end{document}